\documentclass{gtart}

\def\ifplaintex{\expandafter\ifx\csname documentclass\endcsname\relax}

\ifplaintex 
\else
\fi

\def\gtm{{\mathsurround=0pt\it $\cal G\mskip-2mu$eometry \&\ 
$\cal T\!\!$opology $\cal M\mskip-1mu$onographs}}    

\def\gtp{{\mathsurround=0pt\it $\cal G\mskip-2mu$eometry \&\ 
$\cal T\!\!$opology $\cal P\!$ublications}}  

\def\recd{{\small Received:\qua\receiveddate\ifx\reviseddate\relax
\else\qquad Revised:\qua\reviseddate\fi\par}} 

\def\volumenumber#1{\def\thevolumenumber{#1}}
\def\volumeyear#1{\def\thevolumeyear{#1}}
\def\volumename#1{\def\thevolumename{#1}}
\def\papernumber#1{\def\thepapernumber{#1}}
\def\pagenumbers#1#2{\def\startpage{#1}\def\finishpage{#2}}
\def\published#1{\def\publishdate{#1}}
\def\received#1{\def\receiveddate{#1}}
\def\revised#1{\def\reviseddate{#1}}
\def\accepted#1{\def\accepteddate{#1}}

\let\\\par
\let\thevolumenumber\relax\let\thepapernumber\relax
\let\thevolumeyear\relax\let\startpage\relax
\let\finishpage\relax\let\publishdate\relax\let\receiveddate\relax
\let\reviseddate\relax\let\accepteddate\relax\let\theasciititle\relax
\let\theasciiauthors\relax
\let\theasciiabstract\relax

\let\theerratum\relax\let\theasciiemail\relax
\let\theshortauthors\relax\let\theshorttitle\relax

\def\startpage{1}\def\finishpage{15}\def\thepapernumber{77}

\volumenumber{2}
\volumename{Proceedings of the Kirbyfest}
\volumeyear{1999}

\long\def\maketitlep{   

\count0=\startpage

\gtm\nl        
{\small Volume \thevolumenumber: \thevolumename\nl 
\ifx\theerratum\relax\else Erratum \erratumnumber\nl\fi
Pages \startpage--\finishpage\nl}

\vglue 0.1truein   

{\parskip=0pt\leftskip 0pt plus 1fil\def\\{\par\smallskip}{\ifplaintex\large
\else\Large\fi\bf\thetitle}\par\medskip}   
\vglue 0.05truein 

{\parskip=0pt\leftskip 0pt plus 1fil\def\\{\par}{\sc\theauthors}
\par\medskip}%
 
\vglue 0.03truein 

{\small\leftskip 25pt\rightskip 25pt{\bf Abstract}\stdspace\theabstract

{\bf AMS Classification}\stdspace\theprimaryclass
\ifx\thesecondaryclass\relax\else; \thesecondaryclass\fi\par
{\bf Keywords}\stdspace \thekeywords\par}\vglue 7pt

}   

\font\phead=cmsl9 scaled 950
\font=cmsl9 scaled 1050
\font\pnum=cmbx10 scaled 913
\font\lnum=cmbx10 
\font\pfoot=cmsl9 scaled 950
\font=cmsl9 scaled 1050
\ifplaintex
\headline{\vbox to 0pt{\vskip -4.5mm\line{\small\phead\ifnum
\count0=\startpage ISSN 1464-8997 (on line)
1464-8989 (printed) \hfill {\pnum\folio}\else\ifodd\count0\def\\{ }%
\ifx\theshorttitle\relax\thetitle\else\theshorttitle\fi\hfill{\pnum\folio}
\else\def\\{ and }{\pnum\folio}\hfill\ifx\theshortauthors\relax\theauthors
\else\theshortauthors\fi\fi\fi}\vss}}
\footline{\vbox to 0pt{\vglue 0mm\line{\small\pfoot\ifnum\count0=\startpage
Published \publishdate:\qua\copyright\ \gtp\hfill\else
\gtm, Volume \thevolumenumber\ (\thevolumeyear)\hfill\fi}\vss
}}
\else
\makeatletter
\def\@oddhead{{\small\lhead\ifnum\count0=\startpage ISSN 1464-8997 (on line)
1464-8989 (printed) \hfill {\lnum\number\count0}\else\ifodd\count0
\def\\{ }\ifx\theshorttitle\relax \thetitle \else\theshorttitle\fi\hfill
{\lnum\number\count0}\else\def\\{ and }{\lnum\number\count0}
\hfill\ifx\theshortauthors\relax 
\theauthors\else\theshortauthors\fi\fi\fi}}\def\@evenhead{@oddhead}
\def\@oddfoot{\small\lfoot\ifnum\count0=\startpage Published \publishdate:\qua\copyright\ \gtp\hfill\else
\gtm, Volume \thevolumenumber\ (\thevolumeyear)\hfill\fi}
\def\@evenfoot{@oddfoot}
\makeatother
\fi

\let\maketitlepage\maketitlep
\let\makeshorttitle\maketitlepage
\let\maketitle\maketitlepage

\newwrite\gtoutfile
\long\gdef\makeheadfile{  
{\def\\{, }\def\s{ }
\immediate\openout\gtoutfile head.xxx
\immediate\write\gtoutfile{To: math@arxiv.org}
\immediate\write\gtoutfile{Subject: put OR rep NNNNN:ppppp}
\immediate\write\gtoutfile{--text follows this line--}
\immediate\write\gtoutfile{Proxy-for: \ifx\theasciiauthors\relax
\theauthors\else\theasciiauthors\fi\s<\ifx\theasciiemail\relax\theemail\else\theasciiemail\fi>}
\immediate\write\gtoutfile{\noexpand\\}
\immediate\write\gtoutfile{Authors: \ifx\theasciiauthors\relax
\theauthors\else\theasciiauthors\fi}
{\def\\{ }\immediate\write\gtoutfile{Title: \ifx\theasciititle\relax
\thetitle\else\theasciititle\fi}}
\immediate\write\gtoutfile{Subj-class: GT or SG, GR etc}
\immediate\write\gtoutfile{MSC-class: \theprimaryclass\ifx\thesecondaryclass\relax\else, \thesecondaryclass\fi}
\immediate\write\gtoutfile{Journal-ref: Geom. Topol. Monogr. \thevolumenumber\s
(\thevolumeyear) \startpage-\finishpage}
\immediate\write\gtoutfile{Comments: Published by Geometry and Topology Monographs at}
\immediate\write\gtoutfile{\s\s\s  http://www.maths.warwick.ac.uk/gt/GTMon\thevolumenumber/paper\thepapernumber.abs.html}
\immediate\write\gtoutfile{\noexpand\\}
\immediate\write\gtoutfile{}
\ifx\theasciiabstract\relax
\immediate\write\gtoutfile{\theabstract}\else
\immediate\write\gtoutfile{\theasciiabstract}\fi
\immediate\write\gtoutfile{}
\immediate\write\gtoutfile{\noexpand\\}
\immediate\write\gtoutfile{}
\immediate\closeout\gtoutfile}}  

\def\maketitlepage{\maketitlep\makeheadfile}
\let\makeshorttitle\maketitlepage
\let\maketitle\maketitlepage

\volumenumber{4}
\volumename{Invariants of knots and 3-manifolds (Kyoto 2001)}
\volumeyear{2002}
\papernumber{8}
\pagenumbers{103}{117}
\received{30 November 2001}
\revised{12 February 2002}
\accepted{22 July 2002}
\published{28 July 2002}

\usepackage{epsf}

\nocolon

\newtheorem{thm}{Theorem}[section]    
\theoremstyle{definition}
\newtheorem{defn}[thm]{Definition}    

\newtheorem{remark}[thm]{Remark}          

\def\figGroupBa{
\begin{figure}[ht!]
\begin{center}
\mbox{\epsfxsize=2.8in \epsfbox{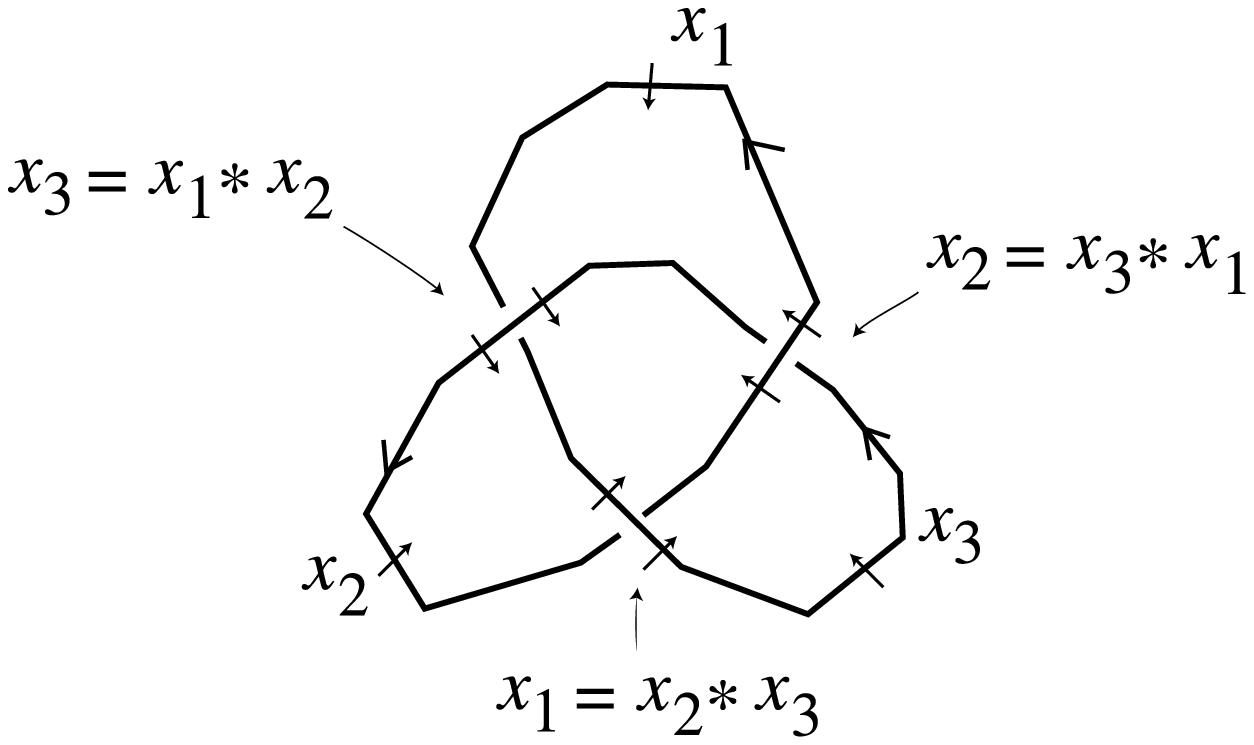}}
\end{center}
\caption{}
\label{fig:GroupBa}\end{figure}} 

\def\figGroupBb{
\begin{figure}[ht!]
\begin{center}
\mbox{\epsfxsize=3in \epsfbox{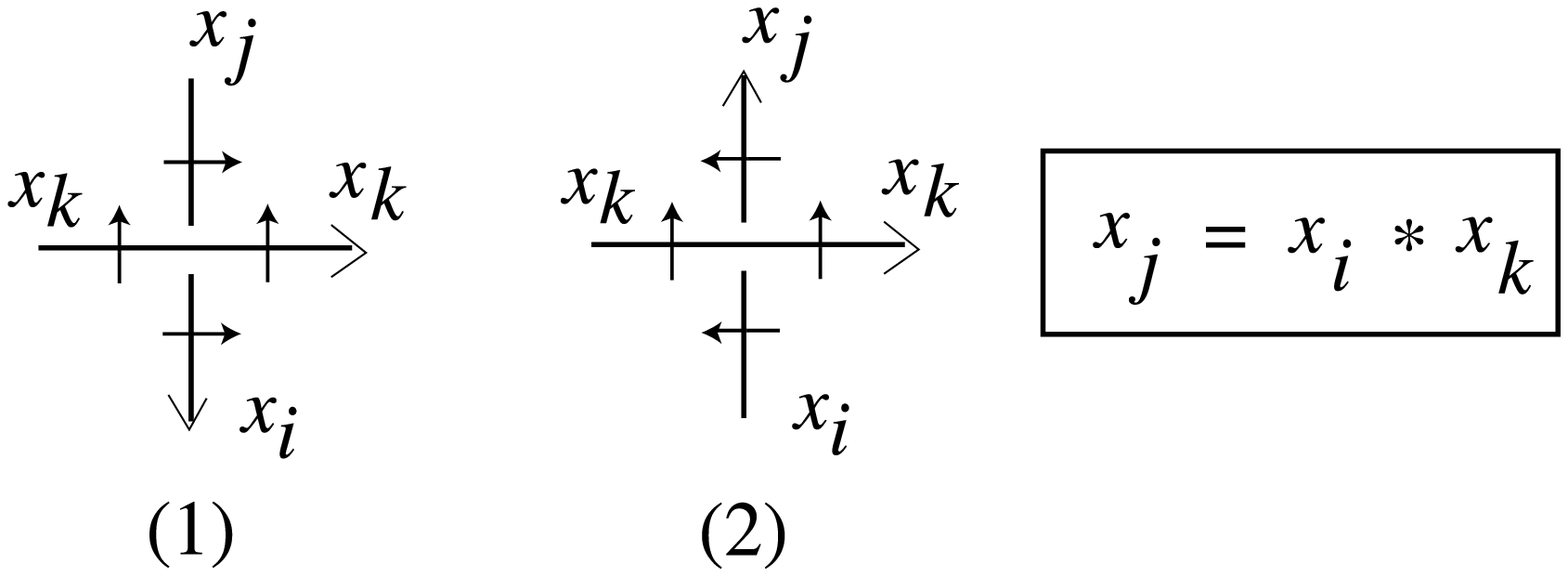}}
\end{center}
\caption{}
\label{fig:GroupBb}\end{figure}} 

\def\figGroupBc{
\begin{figure}[ht!]
\begin{center}
\mbox{\epsfxsize=3in \epsfbox{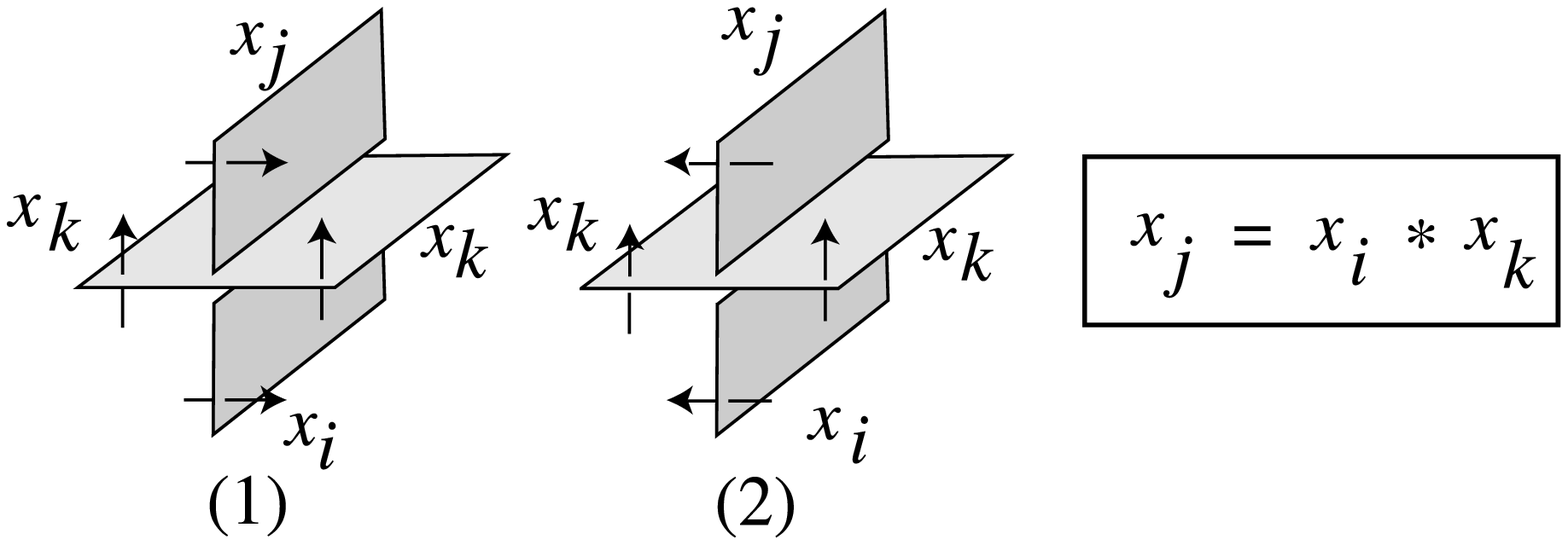}}
\end{center}
\caption{}
\label{fig:GroupBc}\end{figure}}

\begin{document}

\title{Knot invariants derived from quandles and racks}                    
\authors{Seiichi Kamada}                  
\address{Dept. of Math., Hiroshima University, Hiroshima 739-8526, JAPAN}                  
\email{kamada@math.sci.hiroshima-u.ac.jp}                     
\url{http://www.math.sci.hiroshima-u.ac.jp/\char'176kamada/index.html}                     
\begin{abstract}   
The homology and cohomology of quandles and racks are used in knot
theory: given a finite quandle and a cocycle, we can construct a knot
invariant.  This is a quick introductory survey to the invariants of
knots derived from quandles and racks.
\end{abstract}

%
%
%
\makeshorttitle  

%

\section{Introduction}

Many knot invariants can be calculated and interpretated via knot
diagrams.  For example, the (fundamental) quandle of a knot and the
rack of a framed knot \cite{FR, Joyce} can be computed
diagrammatically.
 
In 1942, Takasaki \cite{Takasaki} introduced the notion of a {\it
kei\/}, that is an algebraic object consisting of a non-empty set and
a binary operation on it.  A kei acts on itself by multiplication
{}from the right side (the right translations).  The notion of a kei
is an abstraction of the notion of the symmetric transformations.  In
particular, the right translations of a kei are involutory ($(x*y)*y=
x$ for any elements $x, y$).  When we drop this condition from the
kei, we obtain the notion of a {\it quandle\/}, that was introduced by
Joyce \cite{Joyce} in 1982.  He associated a quandle to a knot, called
the knot quandle and proved that the knot quandle is a complete
invariant of a knot (up to weak equivalence).  At the same time,
Matveev \cite{Matveev} proved a similar result independently.  He
called the notion a {\it distributive groupoid\/}.  Kauffman studied
knot quandles and knot {\it crystals} in \cite{K&P}.  In 1988,
Brieskorn \cite{Brieskorn} introduced the notion of an {\it
automorphic set\/} and, in 1992, Fenn and Rourke \cite{FR} introduced
the notion of a {\it rack\/}.  Automorphic sets and racks are the same
notion, but their actions are the left translations and the right
translations, respectively.  They removed the idempotency condition
($x * x= x$ for any element $x$) {}from the quandle.  Fenn and Rourke
also generalized the notion of the knot quandle to the rack of a
framed knot or a knot diagram.

Although  
knot quandles and racks are strong invariants,
it is not easy to
use  them to distinguish knots by direct calculations.
The situation is similar
to the knot group.
An easy method to use them 
in order to distinguish knots is to
calculate the representations in a given finite
quandle/rack.
Fox's $3$-colorings  are representations of the
knot quandle to the dihedral kei of order $3$ and 
generalized Fox's colorings ($n$-colorings) are
representations to the dihedral kei of order $n$ \cite{CF, FoxTrip, P}.

In the 1990's and the early 2000's,
the homology and cohomology theory for quandles
and  racks appeared \cite{CJKLS, CJKS2, FRS1, FRS2}.
In order to distinguish two quandles or racks,
we may calculate and compare their homology groups.
  Such methods cannot distinguish
a knot and its dual knot (the mirror image with the reversed
orientation), 
 since they have isomorphic knot quandles.
  Carter, Jelsovsky, Kamada, Langford, and Saito
(CJKLS)
\cite{CJKLS}  introduced an invariant of a knot derived from a quandle
cocycle.  Associated with a finite quandle and
a cocycle, we can construct a knot invariant.
The cohomology theory of a quandle introduced
in   \cite{CJKLS} is related to the homology theory of a rack
introdued by Fenn, Rourke and Sanderson (FRS) in
\cite{FRS1,FRS2}, cf \cite{CJKS2, LN} or section \ref{sect:homology}.

The invariant of CJKLS corresponds to the evaluation, by a
fixed cocycle, of the fundamental class reprersented by
a knot diagram in the homology group in the sense of
FRS.   The invariants due to CJKLS and FRS are not
invariants that can be derived only from the knot quandle/rack.
In fact, a trefoil knot and its mirror image have isomorphic quandles
and the $2$-twist spun  trefoil and its orientation reversed
$2$-knot have isomorphic quandles.
However they can be distinguished by
the invariants.

This is a quick introductory survey on the invariants of knots
derived from quandles and racks. The 
author
thanks the organizers of the conference for giving him an oppotunity to give a talk.
The talk was introductory and
concentrated  on the classical knot case according to the organizer's
request.  However this article includes 
some information on the  
$2$-knot case.

\section{Keis, quandles and racks}

A {\it kei}, $X$, is a non-empty set with a binary operation
$(a, b) \mapsto a * b$
satisfying the following axioms: 
\begin{itemize}
\item[(K1)] For any $a \in X$, $a* a =a$.
\item[(K2)] For any $a,b \in X$, $(a*b) *b =a$.
\item[(K3)]
For any $a,b,c \in X$,  $ (a*b)*c=(a*c)*(b*c). $
\end{itemize}

A {\it quandle}, $X$, is a non-empty set with a binary operation
$(a, b) \mapsto a * b$
satisfying the following conditions:
\begin{itemize}
\item[(Q1)] For any $a \in X$, $a* a =a$.
\item[(Q2)] For any $a,b \in X$, there is a unique $c \in X$ such that
$a= c*b$.
\item[(Q3)]
For any $a,b,c \in X$, we have
$ (a*b)*c=(a*c)*(b*c). $
\end{itemize}

The three axioms correspond respectively to the
Reidemeister moves of type I, II, and III
(see \cite{FR}, \cite{K&P}, for example).

A {\it rack}, $X$, is a non-empty set with a binary operation
$(a, b) \mapsto a * b$
satisfying the following conditions:
\begin{itemize} 
\item[(R1)] 
For any $a,b \in X$, there is a unique $c \in X$ such that
$a= c*b$.
\item[(R2)]
For any $a,b,c \in X$, we have
$ (a*b)*c=(a*c)*(b*c). $
\end{itemize}
By definition, a kei is a quandle and
a quandle is a rack:
$$
{\rm \{ keis\} \subset \{quandles \} \subset \{racks\}.  }
$$
Here are some typical examples.

\begin{itemize}
\item
Any set $X$ with the operation $x*y=x$ for any $x,y \in X$ is
a kei, which is called the {\it trivial} kei,
the {\it trivial} quandle or the {\it trivial} rack.
\item
Let $n$ be a positive integer.
For elements  $i, j \in \{ 0, 1, \ldots , n-1 \}$, define
$i\ast j \equiv 2j-i \pmod{n}$.
Then $\ast$ defines a kei
structure  called the {\it dihedral kei},
$R_n$.
This set can be identified with  the
set of reflections of a regular $n$-gon with conjugation
as the kei operation.

\item
A group $X=G$ with
$n$-fold conjugation
as the quandle operation: $a*b=b^{-n} a b^n$.

\item
A $\Lambda (={\bf Z}[T, T^{-1}])$-module $M$
is a quandle with
$a*b=Ta+(1-T)b$, $a,b \in M$, called an {\it  Alexander  quandle}.
For example, $R_4 \cong {\bf Z}_2[T, T^{-1}]/(T^2 +1)$. 
\end{itemize}

See \cite{Brieskorn,FR,Joyce,Matveev}
   for further examples.

A subset $S$ of a rack $X$ is called a {\it
generating  set} if any element of $X$ is obtained from the
elements of $S$ by applying the operation suitably.
If $X$ is a kei, we have
\begin{eqnarray}
\label{eqn:modifykei}
  x * (y * z)
&=& ( (x * (y * z)) * z ) * z \quad\quad {\rm (by \quad K2)}
\nonumber \\
&=& ( (x * z)  * ((y * z) * z ) ) * z \quad\quad {\rm (by \quad K3)}
\nonumber \\
&=& ( (x * z)  * y ) * z \quad\quad {\rm (by \quad K2)}
\end{eqnarray}
for any $x, y, z$ of $X$.  Using this, 
Takasaki \cite{Takasaki} 
and Winker \cite{Winker} (cf.\ \cite{K&P}, p. 195) proved that
any element of a kei $X$ which is expressed by
elements of $S$ can be re-expressed by using the same elements
in the form
as
$$
( \cdots  ( ( (x_1 * x_2) *  x_3)  * x_4 )
\cdots). 
$$
Such an expression is denoted by
$ 
x_1 x_2 x_3 x_4  \cdots $  or $
x_1^{ x_2 x_3 x_4 \cdots}
$ 
following \cite{Takasaki} and \cite{FR}.  
For example, 
\begin{eqnarray}
( a * ( b * c )) * ( d * e)
&=&  ( ( (a * c) *  b)  * c )    * ( d * e) \nonumber \\
&=&  ( ( ( ( (a * c) *  b)  * c ) * e)  *  d)  * e \nonumber \\
&=& a^{cbcede} 
\end{eqnarray}
Fenn and Rourke \cite{FR} introduced this
notation for any rack $X$.
They denote by $a^b$ the element
$a*b$ of $X$ and assume that $a^{b^c}$ stands for
$a^{(b^c)} = a*(b*c)$ and $a^{bc}$ stands for
$(a^b)^c= (a*b)*c$.
Then the axioms are stated as follows:
\begin{itemize}
\item[(R1)] For any $a,b \in X$, there is a unique $c \in X$ such that
$a= c^b$.
\item[(R2)]
For any $a,b,c \in X$, we have
\begin{equation}
\label{eqn:rackid1}
a^{bc}= a^{c b^c}.
\end{equation}
\end{itemize}

The axiom  R1  implies that
the function $S_x: X \to X$ (acting from the right side) defined by
$(u)S_x = u*x = u^x$ is a bijection.
The axiom  R2  implies that it is
a homomorphism; $(u * v)S_x= (u)S_x * (v)S_x$ or
$(u * v)^x= (u^x) * (v^x)$.
We call it 
the {\it right translation} by the element $x$.
    For a word
$  W = x_1^{\epsilon_1}\cdots x_k^{\epsilon_k} $
consisting of elements of $X$ and for an element
$u \in X$, we denote by
$u^W$ the element
$(u)S_{x_1}^{\epsilon_1}\cdots  S_{x_k}^{\epsilon_k}$.
In particular, $a^{b^{-1}}$ stands for the element $c$ with
$c^b =a$.  Then
\begin{equation}
\label{eqn:rackid2}
x^{y^z} = x^{z^{-1}yz}
\end{equation}
for any $x, y, z$ of $X$.
Using this, we see that if $S$ is a generating
set  of a rack, then any element of the rack is expressed by
$u^W$ for some $u \in S$ and a word $W$ consisting of
elements of $S$.  Refer to \cite{FR} for the details.

\section{Knot quandles and colored knot diagrams}

The quandle $Q({\cal K})$ of a knot ${\cal K}$ and 
the rack $R({\cal K})$ of a
framed knot
${\cal K}$ were introduced in \cite{Joyce} and \cite{FR}.
They can be calculated from a diagram.

Let $D$ be a diagram in ${\bf R}^2$ of an oriented knot ${\cal K}$,
and  let $E
=\{x_1, \dots, x_m\}$ be the set of
the arcs of the diagram.  We give each arc of the diagram $D$ a
specific normal direction (co-orientation)
determined  by use of the orientations of
${\cal K}$ and  ${\bf R}^2$.
For a crossing point $\tau$, let $x_i$, $x_j$ and $x_k$ be
the three arcs around $\tau$ such that $x_k$ is
the over-arc at the crossing and $x_i$ is one of the
under-arcs 
away from which the normal direction of $x_k$ points.  
We consider a relation
\begin{equation}
\label{eqn:crossingrelation}
x_j= x_i * x_k
\end{equation}
and call it the {\it crossing relation} at $\tau$ 
(see Fig.~\ref{fig:GroupBb}).

\figGroupBb
 
The kei $K(D)$, the quandle $Q(D)$, or the rack $R(D)$
of the diagram is defined to be the kei, the quandle or the rack
that is generated by the set $E
=\{x_1, \dots, x_m\}$ and the defining relations are
the crossing relations of the crossings of $D$.

The kei $K(D)$ can be defined even if the knot ${\cal K}$ is
unoriented, since 
$x_j= x_i * x_k$ and $x_i= x_j * x_k$ are equivalent
by the second axiom, K2.

By using Reidemeister moves, we see that
$K(D)$ and $Q(D)$ are invariants of the knot ${\cal K}$.
In fact, the quandle $Q(D)$ is
isomorphic to the knot quandle $Q({\cal K})$ of the knot ${\cal K}$, see
\cite{FR, Joyce}.
    Similarly,
$R(D)$ is preserved by Reidemeister moves of type 2 and type 3,
and hence it is an invariant of the framed knot ${\cal K}$,
where we assume
the framing is the blackboard framing of the diagram.
This is
isomorphic to the rack $R({\cal K})$ of the framed knot ${\cal K}$, 
see \cite{FR}.

In general it is difficult to distinguish the 
isomorphism types of two given presentations 
of keis, quandles or racks.  A convenient method 
is to use representations to a finite
kei, quandle or rack $X$.

Let $\rho: K(D) \to X$ be
a homomorphism.  We call it a {\it coloring} of $D$ by $X$.
We color the arc $x_i$ of the diagram $D$
with the element $\rho(x_i)$ of $X$.
The crossing relation $x_j= x_i * x_k$
implies $\rho(x_j)= \rho(x_i) * \rho(x_k)$ in $X$.
This is called the {\it crossing relation} for the coloring
or the {\it coloring condition}.

For example, consider the case that $X$ is the dihedral kei
$R_3 =\{0,1,2\}$ of order three.  The crossing relation
implies that
the three colors $\rho(x_i)$, $\rho(x_j)$ and $\rho(x_k)$
are the same or all of them are distinct.
This is the $3$-coloring condition due to Fox, cf.\ \cite{CF, P}. 
(In Exercise~6 of Chapter VI of \cite{CF}, there is 
an additional condition that all three colors are actually used. 
The $3$-coloring condition together this condition 
is also referred to as Fox's $3$-coloring condition. 
In this paper, we do not assume it.) 

\figGroupBa

A typical diagram of the trefoil knot yields the
presentation
\begin{equation}
\label{eqn:trefoil3gen}
\bigl\langle x_1,   x_2,   x_3 \bigm|
x_3 = x_1 * x_2, \quad
x_1 = x_2 * x_3, \quad
x_2 = x_3 * x_1 \bigr\rangle 
\end{equation} 
(see Fig.~\ref{fig:GroupBa}). 
Then the diagram has $9$ Fox's $3$-colorings.
Three of them are trivial colorings and
the other six are non-trivial colorings.
Since the obvious diagram of the trivial knot has
trivial colorings only, we see that the trefoil knot is not
the trivial knot.

When we use the dihedral kei $R_n$ of order $n$,
we have  Fox's $n$-coloring.  Fox's $n$-coloring is studied in \cite{P}.
In general, by use of a finite kei, we can consider
colorings for unoriented knot diagrams.
When we use a finite quandle or rack $X$ which is not
involutory,
diagrams must be oriented.
A generalized $n$-coloring, called an $(n,r)$-coloring,  was
introduced by
Silver and Williams \cite{SW}. It is
interpretated as a coloring by a certain quandle $X$.

The knot group (the fundamental group of the knot complement)
$G({\cal K})$ can be calculated as the {\it group}
$G(D)$ of a diagram $D$ that is
the group with the generating set $E=\{x_1, \dots, x_m\}$
and the Wirtinger relations
$x_j = x_k^{-1} x_i x_k$ derived from the
crossing points.
The presentation of $G(D)$ is obtained from the
presentation of $Q(D)$ by assuming
$x_j = x_i * x_k$ to be  $x_j = x_k^{-1} x_i x_k$.
In fact, the knot group is a consequence of the knot quandle 
(see below or \cite{FR}).

For a quandle/rack $X$,
assuming $x*y$ to be $y^{-1}xy$, we obtain a group.
It is called the {\it associated group}
of $X$ and denoted by ${\rm As}(X)$. 
More precisely ${\rm As}(X)$ is the group $F(X)/N$ 
where $F(X)$ is the free group generated by the elements of $X$ 
and $N$ is the normal subgroup of $F(X)$ generated by 
the words $(x*y)y^{-1}xy$ for $x, y \in X$, see \cite{FR}. 
A (group) presentation of ${\rm As}(X)$ is obtained from a 
(quandle/rack) presentation of $X$ by reading $x*y$ 
as $y^{-1}xy$, cf.\ Lemma~4.3 of \cite{FR}. 
Conversely, for a group $G$, define a binary operation
$x*y$ by $y^{-1}xy$ and obtain a quandle.
Such a quandle is called the {\it conjugation quandle} of
$G$ and denoted by ${\rm conj}(G)$.

Joyce \cite{Joyce} and Matveev \cite{Matveev} proved the
following theorem,  which shows that the knot
quandle is truely stronger than the knot group
for some composite knots.

\begin{thm}[Joyce \cite{Joyce},  Matveev \cite{Matveev}]
Knots ${\cal K}$
and ${\cal K}'$ with $Q({\cal K}) {\cong} Q({\cal K}')$ are weakly equivalent;
namely, $({\bf R}^3, {\cal K})$ is homeomorphic to
$({\bf R}^3, {\cal K}')$ when we ignore the orientations 
of
${\bf R}^3$ and the knots.
\end{thm}

For an oriented  knot ${\cal K}$, let $-{\cal K}$ denote 
the knot with
the reversed orientation and let  ${\cal K}^*$ denote the mirror image.  
Then $Q({\cal K}) \cong Q(-{\cal K}^*)$.
If the knot ${\cal K}$ is invertible (i.e., ${\cal K} \cong -{\cal K}$), or
amphicheiral (i.e., ${\cal K} \cong {\cal K}^*$), then
$Q({\cal K})$, $Q(-{\cal K})$,  $Q({\cal K}^*)$ and $Q(-{\cal K}^*)$ 
are isomorphic to each other. 
Therefore any invariant which is derived only from the
quandle $Q({\cal K})$ cannot distinguish the trefoil and its mirror image,
since the trefoil is invertible.
The invariants due to CJKLS and FRS can distinguish this pair and hence
they are not invariants that are derived only from the knot quanadle.

\section{Homologies and cohomologies of a quandle and a rack}
\label{sect:homology}

For a rack $X$, 
let $C_n^{\rm R}(X)$ be the free
abelian group generated by
$n$-tuples $(x_1, \dots, x_n)$ of elements of $X$ when
$n$ is a positive integer and put $C_n^{\rm R}(X)=0$.
Define
a homomorphism
$\partial_{n}: C_{n}^{\rm R}(X) \to C_{n-1}^{\rm R}(X)$ by
\begin{eqnarray}
\partial_{n}(x_1, x_2, \dots, x_n)
&=&
\sum_{i=1}^{n} (-1)^{i}\Bigl[
  (x_1, x_2, \dots, x_{i-1}, x_{i+1},\dots,
x_n) \Bigr.
\nonumber \\*
& &
\phantom{\sum_{i=1}^{n}}  
- \Bigl. (x_1^{x_i}, x_2^{x_i}, \dots, x_{i-1}^{x_i}, x_{i+1},
\dots, x_n) \Bigr]
\end{eqnarray}
for $n \geq 2$
and $\partial_n=0$ for
$n \leq 1$.
   Then
$C_\ast^{\rm R}(X)
= \{C_n^{\rm R}(X), \partial_n \}$ is a chain complex.
   Let $C_n^{\rm D}(X)$ be the subset of $C_n^{\rm R}(X)$ generated
by $n$-tuples $(x_1, \dots, x_n)$
with $x_{i}=x_{i+1}$ for some $i \in \{1, \dots,n-1\}$ if $n \geq 2$;
otherwise let $C_n^{\rm D}(X)=0$.  If $X$ is a quandle, then
$\partial_n(C_n^{\rm D}(X)) \subset C_{n-1}^{\rm D}(X)$ and
$C_\ast^{\rm D}(X) = \{ C_n^{\rm D}(X), \partial_n \}$ is a sub-complex of
$C_\ast^{\rm
R}(X)$. 
In this case, let $C_\ast^{\rm Q}(X)$ be the quotient 
complex 
$C_\ast^{\rm R}(X)/ C_\ast^{\rm D}(X)$.
Let $G$ be an abelian group and let ${\rm W}$ stand for  
${\rm R}$, ${\rm D}$ or ${\rm Q}$.   Then we consider the 
homology and cohomology groups   
\begin{eqnarray}
H_n^{\rm W}(X;G)  
= H_{n}(C_\ast^{\rm W}(X) \otimes G), \quad
H^n_{\rm W}(X;G) 
= H^{n}({\rm Hom}(C_\ast^{\rm W}(X), G)). 
\end{eqnarray}

\begin{defn}  
The $n$\/th {\it rack homology} and  {\it rack
cohomology groups\/} \cite{FRS1}
of a rack/quandle $X$ with coefficient group $G$ are 
$
H_n^{\rm R}(X;G) 
$ and $
H^n_{\rm R}(X;G)$. 
The $n$\/th {\it degeneration homology} and  
{\it degeneration cohomology groups\/} of a quandle $X$
with coefficient group $G$ are
$
H_n^{\rm D}(X;G) 
$ and $
H^n_{\rm D}(X;G)$. 
The $n$\/th {\it quandle homology\/} and  
{\it quandle cohomology groups\/} \cite{CJKLS}
of a quandle $X$
with coefficient group $G$ are
$
H_n^{\rm Q}(X;G) 
$ and $
H^n_{\rm Q}(X;G) 
$. 
\end{defn}

The $n$\/th cocycle group of the cochain complex 
${\rm Hom}(C_\ast^{\rm Q}(X), G)$ is called  
the quandle $n$-cocycle group 
and denoted by $Z_n^{\rm Q}(X;G)$.
We will omit the coefficient group $G$  as usual if $G = {\bf Z}$.
   Refer to
\cite{FRS1, FRS2, Flower, Greene}
for some calculations and
applications of the rack homology groups, and
\cite{CJKLS, CJKS1, CJKS2, SSS2, LN, Mochizuki}
for the quandle homology groups.

\begin{thm}[Universal Coefficient Theorem \cite{CJKS2}]
\label{sect:UnivCoefThm}
There exist split exact sequences
\begin{eqnarray}
0 \to H_n^{\rm W}(X) \otimes G \to H_n^{\rm W}(X; G) \to
{\rm Tor}(H_{n-1}^{\rm W}(X), G) \to 0 \\
0 \to {\rm Ext}(H_{n-1}^{\rm W}(X), G) \to H^n_{\rm W}(X; G) \to
{\rm Hom}(H_n^{\rm W}(X), G) \to 0.
\end{eqnarray}
\end{thm}

By definition, there is a natural short exact sequence
\begin{eqnarray}
\label{eqn:basicshortexact}
0 \to C_*^{\rm D}(X) {\to} C_*^{\rm R}(X) {\to}
C_*^{\rm Q}(X) \to 0,
\end{eqnarray}
which is split (in the weak sense), namely,
for each $n$, the short exact sequence  
\begin{eqnarray}
\label{eqn:basicshortexact2}
0 \to C_n^{\rm D}(X) {\to} C_n^{\rm R}(X) {\to}
C_n^{\rm Q}(X) \to 0
\end{eqnarray}
is split.  Therefore we have
a long exact sequence  
\begin{eqnarray}
\label{eqn:basiclongseq}
\cdots \stackrel{\partial_\ast}{\to} H_n^{\rm D}(X;G)
{\to}
H_n^{\rm R}(X;G)
{\to} H_n^{\rm Q}(X;G)
\stackrel{\partial_\ast}{\to} H_{n-1}^{\rm D}(X;G)
\to \cdots.
\end{eqnarray}
In \cite{CJKS2}, it was proved that this sequence 
is split into short exact sequences for small $n$'s and it was
conjectured that so is for all $n$'s.
This was proved by Litherland and Nelson \cite{LN}.

\begin{thm}[Litherland and Nelson \cite{LN}]
The short exact sequence
{\rm (\ref{eqn:basicshortexact})}
is split in the strong sense, that is,
there exist splitting homomorphisms compatible with the
boundary maps.  Thus, for each $n$,  there is a short exact sequence
\begin{eqnarray}
0 {\to} H_n^{\rm D}(X;G)
{\to}
H_n^{\rm R}(X;G)
{\to} H_n^{\rm Q}(X;G)
  {\to} 0. 
  \end{eqnarray}
\end{thm}

Suppose that $X$ is a quandle with finitely many orbits, where
orbits means orbits
by the inner-automorphisms of $X$, i.e.,
orbits by the right action of $X$
on itself by the right translations. (In \cite{Joyce}
orbits by 
all the  
automorphisms of $X$ are studied.)
Let ${\bf Z}[{\rm Orb}(X)]$ be  
the free abelian group generated by the orbit set, ${\rm Orb}(X)$,
of $X$.  In \cite{CJKS2, Greene} it is shown that
$$ H_1^{\rm D}(X) =0, \quad
H_1^{\rm R}(X) \cong H_1^{\rm Q}(X)
\cong
{\bf Z}[{\rm Orb}(X)],  \quad {\rm and} \quad
H_2^{\rm D}(X) \cong
{\bf Z}[{\rm Orb}(X)]. $$
Hence
$$H_2^{\rm R}(X) \cong H_2^{\rm Q}(X) \oplus
{\bf Z}[{\rm Orb}(X)]. $$
In \cite{LN}(Theorem~2.2)  it is proved that
$$H_3^{\rm R}(X) \cong H_3^{\rm Q}(X) \oplus H_2^{\rm Q}(X)
\oplus
{\bf Z}[{\rm Orb}(X) \times {\rm Orb}(X)]. $$
The quandle (co)homology groups of some finite 
quandles were caluculated in \cite{CJKS1} (Table 1).
However there are errata due to the authors.  They thank  
R. A. Litherland, S. Nelson \cite{LN} and
T. Mochizuki \cite{Mochizuki}
for pointing out and correcting them. 
Further calculations on homology/cohomology groups of 
finite Alexander quandles are
found in their papers \cite{LN} and \cite{Mochizuki}.

\section{Knot invariants}

Let $X$ be a finite quandle, let $G$ be an abelian group
(written in multiplicative notation), and let
$\phi \in C^2_{\rm Q}(X; G)$ be a quandle $2$-cocycle.
 
Let $D$ be an oriented knot diagram.  Recall that
a homomorphism $\rho: Q(D) \to X$ is called a coloring and
satisfies the crossing condition
$\rho(x_j)= \rho(x_i) * \rho(x_k)$
at each crossing, where $x_i$, $x_j$ and $x_k$ are as before 
(see Fig.~\ref{fig:GroupBb}).
For each crossing point $\tau$, we consider an element
$W_\phi(\tau, \rho)$ of $G$ determined by
\begin{equation}
W_\phi(\tau, \rho) = \phi(\rho(x_i), \rho(x_k))^\epsilon,
\end{equation}
where $x_i$ and $x_k$ are the arcs around $\tau$
corresponding the crossing relation $x_j = x_i * x_k$,
and $\epsilon$ is the sign of the crossing point $\tau$.
Consider the element
\begin{equation}
\Phi_\phi(D) = \sum_{\rho} \prod_{\tau} W_\phi(\tau, \rho)
\end{equation}
of the group ring ${\bf Z}G$, where $\rho$ runs all colorings
of $D$ by $X$ and $\tau$ runs all crossing points of $D$.
   Since $\phi$ is a cocycle,  the value $\Phi_\phi(D)$ is
preserved by Reidemeister moves.  Therefore it is an
invariant of the knot ${\cal K}$ represented by the diagram.
Thus we denote this value by $\Phi_\phi({\cal K})$ and call it
the {\it state-sum invariant\/} or the
{\it quandle cocycle invariant\/}, \cite{CJKLS}.
If $\phi$ and $\phi'$ are cocycles that are 
cohomologous 
then
$\Phi_\phi({\cal K})= \Phi_{\phi'}({\cal K})$.  Thus this invariant
depends on the quandle cohomology class of $\phi$.
    When $\phi$ is a trivial cocycle, the element $W_\phi(\tau, \rho)$
is the trivial element of $G$.  In this case the invariant
$\Phi_\phi({\cal K})$ is equal to the number of all colorings of
$D$ by $X$.
In \cite{SSSxx, CENS}, the {\it abelian extension} $E =E(X,G, \phi)$ is 
defined as $G \times X$ as a set, with the quandle operation 
$(g_1, x_1)*(g_2, x_2) = (g_1 \phi(x_1,x_2), x_1*x_2)$. 
The state-sum $\Phi_\phi(D)$ is an obstruction to extending 
the colorings of $D$ by $X$ to colorings by $E(X,G, \phi)$, 
see Theorems~4.1 and 4.5 of \cite{CENS}.

Let $X$ be a finite quandle, let $G$ be an abelian group
(written in multiplicative notation), and let
$\theta \in C^3_{\rm Q}(X; G)$ be a quandle $3$-cocycle.

A {\it shadow coloring} (or {\it face coloring}) of
$D$ by $X$ is a
function $\tilde{\rho} : \tilde{E} \rightarrow X$, where
$ \tilde{E}$ is the set of arcs of the diagram
$D$ and the regions separated
by the underlying immersed curve of the diagram,
satisfying the coloring condition
($\tilde{\rho}(x_j)= \tilde{\rho}(x_i) * \tilde{\rho}(x_k)$) and the
condition that
if $y_i$ and $y_j$ are regions  
which are adjacent to
an arc $x_k$ and the normal direction of $x_k$ points away
{}from $y_i$ toward $y_j$, then
$$ \tilde{\rho}(y_j) = \tilde{\rho}(y_i) * \tilde{\rho}(x_k). $$
We call this condition the {\it face coloring condition}.

Let $\tau$ be a crossing point and let $x_i$ and $x_k$ be
the arcs around $\tau$ appearing in the
crossing relation $x_j = x_i * x_k$  at $\tau$ as before.
Let $y$ be the region which is one of the four regions
around $\tau$ such that the normal directions of
$x_i$ and $x_k$ are away from the region
$y$.
Then we give this crossing point $\tau$  an element
${W}_\theta(\tau, \tilde{\rho})$ of $G$ determined by
\begin{equation}
{W}_\theta(\tau, \tilde{\rho})
= \theta(\tilde{\rho}(y), \tilde{\rho}(x_i),
\tilde{\rho}(x_k))^\epsilon,
\end{equation}
where $\epsilon$ is the sign of the crossing point.
Consider  
the element
\begin{equation}
\Psi_\theta(D) = \sum_{\tilde{\rho}} \prod_{\tau}
{W}_\theta(\tau, \tilde{\rho})
\end{equation}
of the group ring ${\bf Z}G$, where $\tilde{\rho}$ runs all
shadow colorings
of $D$ by $X$ and $\tau$ runs all crossing points of $D$.
   Since $\theta$ is a cocycle, the value $\Psi_\theta(D)$ is
preserved by Reidemeister moves.  It is an
invariant of the knot ${\cal K}$ represented by $D$.
We denote this value by $\Psi_\theta({\cal K})$.
If $\theta$ and $\theta'$ are cohomologous, then
$\Psi_\theta({\cal K})= \Psi_{\theta'}({\cal K})$.  Thus this invariant
depends on the quandle $3$-cohomology class of $\theta$.
    When $\theta$ is a trivial cocycle, the invariant
$\Psi_\theta({\cal K})$ is equal to the number of the shadow colorings of
$D$ by $X$.

Colorings and shadow colorings are defined for
oriented knotted surfaces
in $4$-space similarly using their diagrams in $3$-space.
For a knotted surface ${\cal K}$  in ${\bf R}^4$, by
modifying it slightly, we may assume that the
projection $p : {\cal K} \to {\bf R}^3$ to the $3$-space is
a generic map.  The singularity set of the projection consists 
of double points, triple points and branch points.
By removing a 
small regular neighborhood of the under-curve
of the double curve, we have a compact embedded surface in
the $3$-space.  This is a diagram of the knotted surface ${\cal K}$
(see \cite{CS:book} for details).  Let $D$ be the diagram and let
  $E= \{x_1, \dots, x_m\}$ be the set of sheets of the
diagram.
Using the orientation of ${\cal K}$ and the orientation of $3$-space,
we give a normal direction (co-orientation) to
each sheet.  At a double curve of the
projection $p({\cal K})$, let  $x_i$, $x_j$ and $x_k$ be the
three sheets around the double curve such that
$x_k$ is the over-sheet and $x_i$ and $x_j$ are under-sheets
and that the normal direction is away from $x_i$ toward $x_j$.
At each double curve, we consider a relation
\begin{equation}
   x_j = x_i * x_k
\end{equation}
and call it the {\it crossing relation} around the double curve 
(see Fig.~\ref{fig:GroupBc}).
The {\it quandle} $Q(D)$ of the diagram is defined to be the
quandle generated by $E= \{x_1, \dots, x_m\}$ and the
defining relations are the crossing relations around the
double curves.  This quandle is isomorphic to the
quandle $Q({\cal K})$ of the knotted surface ${\cal K}$ in the sense of
\cite{FR, Joyce}.

\figGroupBc

A homomorphism
$\rho: Q(D) \to X$ is called a {\it coloring} of
the diagram $D$ by $X$.  This is equivalent to coloring
the sheets by elements of $X$ such that the
colors $\rho(x_i)$, $\rho(x_j)$ and $\rho(x_k)$
of the three sheets $x_i, x_j$ and $x_k$ around a double
curve satisfies
$ \rho(x_j) = \rho(x_i) * \rho(x_k). $

Let $X$ be a finite quandle, let $G$ be an abelian group
(written in multiplicative notation), and let
$\theta \in C^3_{\rm Q}(X; G)$ be a quandle $3$-cocycle.
   For each triple point $\tau$ of the diagram $D$,
we consider an element
$W_\theta(\tau, \rho)$ of $G$ determined by
\begin{equation}
W_\theta(\tau, \rho) = \phi(\rho(x_i), \rho(x_j), \rho(x_k))^\epsilon,
\end{equation}
where $x_i$, $x_j$ and $x_k$ are the sheets around $\tau$
such that $x_k$ is the upper-sheet, $x_j$ is one of the two
middle-sheets from which the normal direction of
$x_k$ is away, and $x_i$ is one of the four lower-sheets away  
from which the normal directions of $x_k$ and $x_j$ point,  
and $\epsilon$ is the sign of the triple point $\tau$.
Consider the element
\begin{equation}
\Phi_\theta(D) = \sum_{\rho} \prod_{\tau} W_\phi(\tau, \rho)
\end{equation}
of the group ring ${\bf Z}G$, where $\rho$ runs all colorings
of $D$ by $X$ and $\tau$ runs all triple points of $D$.
  Since $\phi$ is a cocycle,  the value $\Phi_\phi(D)$ is
preserved by Roseman moves (\cite{CS:book, Rose}).  Therefore it is an
invariant of the knotted surface ${\cal K}$ represented by the diagram.
We denote this value by $\Phi_\phi({\cal K})$.
If $\phi$ and $\phi'$ are cohomologous, then
$\Phi_\phi({\cal K})= \Phi_{\phi'}({\cal K})$.  Thus this invariant
depends on the quandle cohomology class of $\phi$.
    When $\phi$ is a trivial cocycle,  the invariant
$\Phi_\phi({\cal K})$ is equal to the number of the colorings of
$D$ by $X$.

\section{Examples}

Let $X$ be the dihedral kei
$R_3=\{0,1,2\}$ of order three, let $G$ be the cyclic group
$\langle t \, | \,  t^3=1 \rangle$, and let
$\theta \in Z^3_{\rm Q} (R_3; {\bf Z}_3)$ be the
quandle $3$-cocycle defined by
$$\theta = \chi_{012} \chi_{021}
\chi_{101} \chi_{201} \chi_{202}\chi_{102}$$
where
$\chi_{abc} (x,y,z) = t$ if $(x,y,z)=(a,b,c)$; 
otherwise $\chi_{abc} (x,y,z) = 1$.  
   In this situation, we consider the invariant 
$\Psi_\theta$ for classical
knots (defined by shadow colorings) and  
the invariant 
  $\Phi_\theta$ for knotted surfaces
(defined by colorings). 
Note that the values 
$\Psi_\theta({\cal K})$ and $\Phi_\theta({\cal K})$ are
elements of
${\bf Z}\langle t \, | \,  t^3=1 \rangle$.

\begin{thm}[Rourke and Sanderson \cite{RS}]
The invariant $\Psi_\theta$ distinguishes the trefoil $T(2,3)$ and
its mirror image $T(2,-3)$;
$\Psi_\theta( T(2, 3) )  \neq \Psi_\theta( T(2, -3) )$.
\end{thm}

\begin{thm}[Carter et al.\ \cite{SSSx}]
A torus knot/link $T(2, n)$ is $3$-colorable if and only if 
$n=3k$ for some integer $k$.  In this case, 
$\Psi_\theta( T(2, 3k) ) = 9 + 18 t^k $.
\end{thm}

\begin{thm}[Carter et al.\ \cite{SSSx}]
A torus knot/link $T(3, n)$ is $3$-colorable if and only if 
$n=2k$ for some integer $k$.  In this case, 
if $k$ is not a multiple of $3$, then
$\Psi_\theta( T(3, 2k) ) = 9 +  18 t^k $; otherwise 
$\Psi_\theta( T(3, 2k) ) =45$.
\end{thm}

There are four knots in the table with less than $8$ crossings that
are $3$-colorable: $3_1$, $6_1$,  $7_4$, and $7_7$.
$\Psi_\theta(3_1)= \Psi_\theta(7_4)= \Psi_\theta(7_7)= 9 + 18t$
and $\Psi_\theta(6_1)= 27$, \cite{SSSx}.
Further examples are calculated in  \cite{CJKS1, SSSx}.

\begin{thm}[Carter et al.\ \cite{CJKLS},
Rourke and Sanderson \cite{RS}]
Let ${\cal K}$ be the\break $2$-twist spun trefoil and
let $-{\cal K}$ be the same $2$-knot with the reversed orientation.
Then
$\Phi_\theta({\cal K}) = 3+6 t^2$ and $\Phi_\theta(-{\cal K}) = 3+6 t$.
Thus the invariant  $\Phi_\theta$ distinguishes
${\cal K}$ and $-{\cal K}$.
\end{thm}

Carter et al.\   \cite{CJKS1} (Theorem 7.2, Cor. 7.3 and 7.4) and
Satoh \cite{Satoh} (Theorem~12) gave formulae on the invariants
$\Phi_\theta({\cal K})$ for the twist-spun knots. The following theorem is
a consequence of the formulae.

\begin{thm}[Carter et al.\ \cite{CJKS1}, Satoh\cite{Satoh}]
Let ${\cal K}$ be the $m$-twist spun trefoil.
Then
\begin{eqnarray}
\Phi_\theta({\cal K}) =
\left\{
\begin{array}{ll}
3+6 t^2  &{\rm for }\quad m\equiv 2 {\rm ~mod~} {6} \\
3+6 t  &{\rm for }\quad m\equiv 4 {\rm ~mod~} {6} \\
9  &{\rm for }\quad m\equiv 0 {\rm ~mod~} {6} \\
3  &{\rm otherwise. }
\end{array}
\right.
\end{eqnarray}
\end{thm}

\begin{remark}
{\rm
(1) Here we used quandle cohomology of a finite quandle
to define the invariants $\Phi$ and $\Psi$.  We may use the
rack cohomology of a quandle or a rack
to define a framed knot.  For the purpose 
of constructing a $3$-manifold invariant, 
one should use the rack homology.

(2)
To define the invariants $\Phi$ and $\Psi$, we used
the summation over all possible colorings. When we need
an invariant of a colored knot, we do not need to take
the summation.
} \end{remark}

\begin{remark}
{\rm
(1) Shadow coloring of a classical knot
diagram  is closely related
to coloring a knotted surface diagram and its lower decker set.
The correspondence is given in \cite{SSS2}.

(2) Recall that a shadow coloring of a classical knot diagram
is a function $\tilde{\rho} : \tilde{E} \rightarrow X$.
The arcs and the regions in $\tilde{E} $ are colored by elements of
$X$.  However, when we use an $X$-set $Y$ in the sense of
Fenn et al. (cf.\ \cite{FRS1, FRS2}), we may color the regions
with elements of $Y$ and the arcs with elements of $X$.
Equivalently, we may introduce and use the cohomologies with twisted
coefficients.  So there are a lot of variations of the invariants
derived from quandle/rack (co)homologies. 
For example, see \cite{CES}. 
  } \end{remark}

\begin{remark}
{\rm Recently, S. Satoh and A. Shima \cite{SatohShima} proved that
if $\Phi_\theta({\cal K})$ is not an integer for a knotted
surface ${\cal K}$, then
the minimum triple point number $t({\cal K})$ of the
generic projections of ${\cal K}$ is greater than three, where
  $\Phi_\theta$ is the invariant used in this section.
As a corollary, they proved that $t({\cal K})=4$
for the $2$-twist spun trefoil ${\cal K}$.
} \end{remark}

%
%

%
%
%
\Addresses
%
%
%
%
\end{document}